\documentclass[11pt]{article}
\usepackage{amsmath,amsfonts,amsthm,amssymb,dsfont,epsfig}
\newtheorem{theorem}{Theorem}[section]
\newtheorem{corollary}[theorem]{Corollary}
\newtheorem{lemma}[theorem]{Lemma}

\newtheorem{definition}[theorem]{Definition}
\newtheorem{claim}[theorem]{Claim}

\newcommand\bb{\boldsymbol{\mathit{b}}}

\newcommand\dd{\boldsymbol{\mathit{d}}}

\newcommand\uu{\boldsymbol{\mathit{u}}}
\newcommand\vv{\boldsymbol{\mathit{v}}}

\newcommand\xx{\boldsymbol{\mathit{x}}}

\newcommand\yy{\boldsymbol{\mathit{y}}}

\newcommand\zz{\boldsymbol{\mathit{z}}}

\def\Reals#1{\mathds{R}^{#1}}

\newcommand{\la}{\langle}
\newcommand{\ra}{\rangle}
\newcommand{\lmx}{\left[\begin{matrix}}
\newcommand{\rmx}{\end{matrix}\right]}
\newcommand{\cleq}{\preccurlyeq}

\newcommand{\mybox}[1]{\medskip \noindent\fbox{\parbox{6.5in}{#1}}\medskip}
\newcommand{\comment}[1]{}

\begin{document}
\title{Support-Graph Preconditioners for 2-Dimensional Trusses}

\author{
Samuel I. Daitch\thanks{
Partially supported by NSF grant CCR-0324914.}\\
Department of Computer Science \\
Yale University \\
\and
Daniel A. Spielman\thanks{
Partially supported by NSF grant CCR-0324914.}\\
Department of Computer Science \\
Yale University \\
}


\maketitle
\thispagestyle{empty}

\begin{abstract}
We use support theory, in particular the fretsaw extensions of
  Shklarski and Toledo~\cite{Fretsaw}, to design preconditioners for
  the stiffness matrices of 
  2-dimensional truss structures that are stiffly connected.
  Provided that all the lengths of the trusses are within constant factors
  of each other, that the angles at the corners of the triangles are bounded
  away from 0 and $\pi$, and that the elastic moduli 
  and cross-sectional areas of all the truss elements
  are within constant factors of each other, our preconditioners allow
  us to solve linear equations in the stiffness matrices to 
  accuracy $\epsilon$ in time
  $O (n^{5/4} (\log^{2}n \log \log n)^{3/4} \log (1/\epsilon ))$.
\end{abstract}

%
%





\section{Preconditioning}

When solving a linear system in an $n\times n$ positive semidefinite matrix
$A$, the running time of an iterative solver can often be sped up 
by {\bf supporting} $A$ with another matrix $B$, called a 
{\bf preconditioner}.  An effective preconditioner $B$ has the properties
that it is much easier to solve than $A$, and that
$A$ has a low condition number relative to $B$.

We define here generalized eigenvalues and condition numbers:

\begin{definition}
For positive semidefinite $A,B$,
the {\bf maximum eigenvalue},
{\bf minimum eigenvalue}, and
{\bf condition number} of $A$ relative to $B$ 
are defined respectively as
$$
\lambda_{max}(A,B) = 
\max_{\xx:\xx\perp null(B)}\frac{\xx^TA\xx}{\xx^TB\xx}
$$
$$
\lambda_{min}(A,B) = 
\min_{\xx:\xx\perp null(A)}\frac{\xx^TA\xx}{\xx^TB\xx}
$$
$$
\kappa(A,B) = \lambda_{max}(A,B)/\lambda_{min}(A,B)
$$
where 
$\xx\perp null(S)$ means that $\xx$ is orthogonal to the null space
of $S$.
\end{definition}
Note that the standard condition number of $A$ can be expressed as
$\kappa(A)=\kappa(A,I)$.

The conjugate gradient method is an example of a linear solver that
can be sped up using a preconditioner.  The precise analysis of the
running time can be found, for example, in \cite{Axelsson}:

\begin{theorem}[\cite{Axelsson}]
For positive semidefinite $A,B$, and vector $\bb$, let $\xx$ satisfy
$A\xx=\bb$.  Each iteration of the preconditioned conjugate gradient method
multiplies one vector by $A$, solves one linear system in $B$, and performs a
constant number of vector additions.  For $\epsilon>0$,
it requires at most $O(\sqrt{\kappa(A,B)}\log(1/\epsilon))$ such
iterations to produce a 
$\tilde{\xx}$ that satisfies
$$
\|\tilde{\xx}-\xx\|_{A}\leq \epsilon\|\xx\|_{A}
$$
\end{theorem}

\subsection{Using a Larger Matrix}

In certain situations it may be easier to find a good preconditioner 
for a matrix $A$ if we
treat $A$ as being larger than it really is.
That is, if we pad $A$
with zeros to form a larger square matrix
$A' = \left[\begin{matrix} A & 0 \\ 0 & 0 \end{matrix}\right]$,
it may be simpler to find a good preconditioner $B$ for $A'$.
We then need to show how to use $B$ to yield a preconditioner 
for the original matrix $A$. To this end, we define the Schur complement:
\begin{definition}
For square matrices 
$A$ and 
$B = 
\left[\begin{matrix} B_{11} & B_{12} \\ B_{12}^T & B_{22} \end{matrix}\right]$,
where square submatrix $B_{11}$ is the same size as $A$,
and such that $B_{22}$ is nonsingular,
the {\bf Schur complement} of $B$ with respect to $A$ is
$$B_S = B_{11} - B_{12}B_{22}^{-1}B_{12}^T$$
\end{definition}

While $B_S$ will not automatically be a good preconditioner for $A$ simply
because $B$ is a good preconditioner for $A'$, we do know that the maximum 
eigenvalue will be the same:

\begin{lemma}\label{lambdamax}
For positive semidefinite $A,B$,
$\lambda_{max}(A',B) = \lambda_{max}(A,B_S)$.
\end{lemma}

We also know that solving a linear system in $B_S$ is as easy as solving a
linear system in $B$:

\begin{lemma}\label{equivlem}
$B\lmx \xx \\ \yy \rmx = 
\lmx \bb \\ 0 \rmx$ implies
$B_S\xx = \bb$
\end{lemma}

For completeness, we give proofs for these lemmas in Appendix \ref{apx:lemmas}.

\subsection{Congestion-Dilation}

Suppose that we have matrices $A$ and $B$ that can be expressed as the sums
of other matrices, i.e. $A=\sum_i A_i$ and $B=\sum_j B_j$, and that
we know how to support each $A_i$ by a subset of the $B_j$ matrices.
In this situation, we can use the following lemma to show how
$B$ supports $A$:

\begin{lemma}[Congestion-Dilation Lemma]\label{cong-dil}

Given the symmetric positive semidefinite matrices
$A_1,...,A_n,B_1,...,B_m$ and $A = \sum_i A_i$ and $B = \sum_j B_j$
and given sets $\Sigma_i \subseteq [1,...,m]$ and real values $s_i$
that satisfy

$$
\lambda_{max}(A_i,\sum_{j\in \Sigma_i} B_j) \leq s_i
$$

it holds that

$$
\lambda_{max}(A,B) \leq \max_j \left(\sum_{i:j\in\Sigma_i} s_i\right)
$$

\end{lemma}

A proof of this lemma is given in Appendix \ref{apx:lemmas}.
It is an adaptation of, for example, Proposition 9.4 in \cite{SupportTheory}.
The expression $\sum_{i:j\in\Sigma_i} s_i$ can be though of as the
{\bf congestion} of matrix $B_j$, analogous to the concept of congestion for
graph embeddings which we define in Section \ref{stretch}.

\section{Trusses and Stiffness Matrices}

\begin{definition}
A {\bf 2-dimensional truss} $\mathcal{T}=\la n,\{\vv_{i} \}_{i=1}^{n},E,\gamma\ra$
is an undirected weighted planar graph
with vertices $[n]=\{1,\dotsc ,n \}$ and edges $E$, 
with vertex $i\in [n]$ embedded at point $\vv_i\in \Reals2$.  
We allow multiple vertices to be 
embedded at the same point.

An edge $e=(i,j)\in E$,
also called a {\bf truss element},
represents a straight idealized bar from $\vv_{i}$ to $\vv_{j}$,
with positive weight $\gamma(e)$ denoting the product of
the bar's cross-sectional area and the elastic modulus of its material.
A {\bf truss face} is a triple $\{i,j,k\}$ such that 
$\{(i,j),(i,k),(j,k) \}\subseteq E$ and
no vertex is in
the interior of the triangle formed by $\vv_{i},\vv_{j},\vv_{k}$.
Every truss element is required to be contained in some truss face.
\end{definition}

There is a particular type of linear system that arises
when analyzing the forces on a truss using the finite element method.
We define here the type of matrix we wish to solve:

\begin{definition}
Given a truss $\mathcal{T}=\la n,\{\vv_{i} \}_{i=1}^{n},E,\gamma\ra$, 
for each truss element $e=(i,j)\in E$ 
we define a length $2n$ column vector $\uu_e=[u_e^1\ ...\ u_e^{2n}]^T$ 
with 4 nonzero entries satisfying
$[ u_e^{2i-1}\  u_e^{2i} ]^T =
-[ u_e^{2j-1}\  u_e^{2j} ]^T =
\frac{\vv_i-\vv_j}{|\vv_i-\vv_j|}$,
and we define the $2n\times 2n$ matrix
\[
A_{e} = \frac{\gamma(e)}{|\vv_i-\vv_j|} \uu_e \uu_e^T
\]
The {\bf stiffness matrix} of the truss is then given by:
\[
A_\mathcal{T} = \sum_{e\in E} A_{e}
\]
\end{definition}

Note that a stiffness matrix is positive semidefinite,
since for all $\xx$ we have
\[
\xx^TA_\mathcal{T}\xx = \sum_{e=(\vv_i,\vv_j)\in E} 
\frac{\gamma(e)}{|\vv_i-\vv_j|} \xx^T\uu_e \uu_e^T\xx
= \sum_{e=(\vv_i,\vv_j)\in E} 
\frac{\gamma(e)}{|\vv_i-\vv_j|} (\xx^T\uu_e)^{2}
\geq 0
\]

We would like to restrict our attention to trusses with a unique, well-behaved
stress-free position.
To this end, we make the following definitions:
\begin{definition}
The {\bf rigidity graph} $Q_\mathcal{T}$
of a truss $\mathcal{T}$ is the graph with vertex set given by the
set of truss faces of $\mathcal{T}$, and with edges connecting faces that
share an edge.

We say that a truss $\mathcal{T}$ is {\bf stiffly-connected} if
(1) $Q_\mathcal{T}$ is connected, and (2)
for every $i\in [n]$, $Q_\mathcal{T}^{i}$ is connected, where
$Q_\mathcal{T}^{i}$ is the graph induced by $Q_\mathcal{T}$ on the set of
faces containing vertex $i$.
\end{definition}

The main contribution of this paper is an algorithm {\tt TrussSolver}
for solving linear systems in stiffness matrices of stiffly-connected trusses.
We will describe the algorithm later, but we state here the result
of our analysis of the running time:

\begin{theorem}[Main Result]\label{mainthm}
For any stiffly-connected truss 
$\la n,\{\vv_{i} \}_{i=1}^{n},E,\gamma\ra$
such that
\begin{itemize}
\item all truss elements have lengths in the range $[l_{min},l_{max}]$
\item all angles of truss faces are in the range $[\theta_{min},\pi-\theta_{min}]$.
\item all weights are in the range $[\gamma_{min},\gamma_{max}]$.
\end{itemize}
for positive constants $l_{min},l_{max},\theta_{min},\gamma_{min},\gamma_{max}$,
{\tt TrussSolver} solves linear systems in matrix $A_\mathcal{T}$
within relative error $\epsilon$
in time
$O\left(n^{5/4}(\log^2 n\log\log n)^{3/4}\log(1/\epsilon)\right)$
\end{theorem}


\subsection{Fretsaw Extension}

We will precondition the stiffness matrix using a fretsaw extension,
a technique described in \cite{Fretsaw}.  The fretsaw extension of a truss
is a new truss created by splitting some of the vertices into multiple copies,
without changing the identity of the truss faces.

\begin{figure}
\begin{center}
\epsfig{file=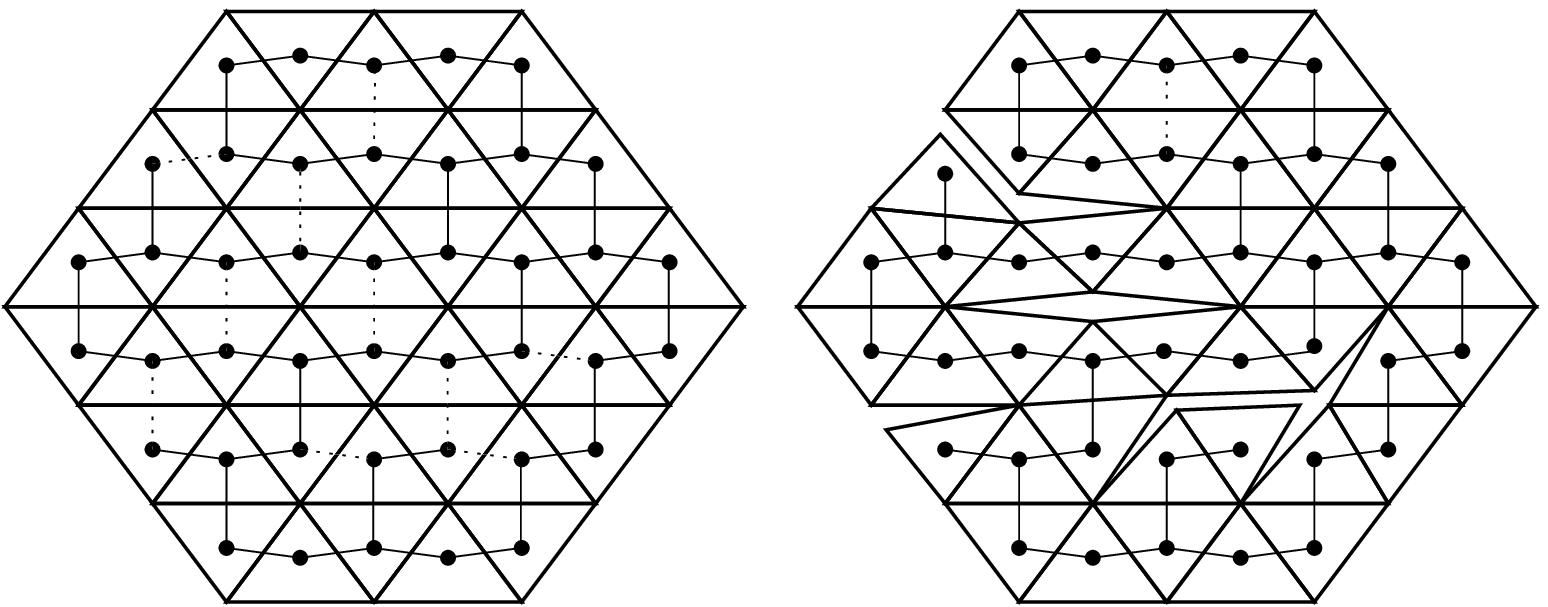, height=1.25in}
\end{center}
\caption{The truss on the right is a fretsaw extension of the truss on the
  left, as given by the {\tt fretsaw} algorithm.
  The vertex positions in the fretsaw extension are distorted slightly so as
  to be able to distinguish vertex copies in the same location.
  The subgraph F is shown as
  solid lines, while the rest of the trusses' connectivity graphs are shown as
  dotted lines.  Note that the connectivity graph of the fretsaw extension has
  one edge not in F.}
\label{figure:truss}
\end{figure}

\begin{definition}
Let $\mathcal{T}=\la n,\{\vv_{i} \}_{i=1}^{n},E,\gamma\ra$ and 
$\mathcal{T}'=\la m,\{\vv'_{i} \}_{i=1}^{m},E',\gamma'\ra$ be 2-dimensional trusses,
with $m>n$.
Let $F$ and $F'$ be the sets of truss faces of $\mathcal{T}$ and $\mathcal{T}'$ respectively,
Let $\rho:F'\rightarrow F$ be a bijection between the faces,
and let $\pi:[m]\rightarrow[n]$ be a surjection on the vertices.
 
$\mathcal{T}'$ is the {\bf $(\rho,\pi)$-fretsaw extension} of $\mathcal{T}$ if
\begin{itemize}
\item for all $i\in [m]$, the $i$th vertex of $\mathcal{T}'$ is a copy of 
the $\pi (i)$th vertex of $\mathcal{T}$, i.e.
$\vv'_{i} = \vv_{\pi(i)}$
\item for all faces $f=(i,j,k)\in F'$,
the vertices of face $f$ are copies of the vertices of $\rho(f)$, i.e.
$\rho (i,j,k) = (\pi(i),\pi(j),\pi(k))$
\item every edge $(i,j)\in E'$ has the same weight as the edge of which it is a copy,
i.e. $\gamma'(i,j)=\gamma(\pi(i),\pi(j))$
\end{itemize}
Since each vertex in $\mathcal{T}$ has at least one copy in $\mathcal{T}'$,
we follow that convention that $\forall i\in[n]: \pi(i) = i$.
Thus, for $i\in [n]$, vertex $i$ in $\mathcal{T}'$ can be considered 
the ``original copy'' of vertex $i$ in $\mathcal{T}$.
\end{definition}
A fretsaw extension has the following property which makes it
a useful preconditioner:

\begin{lemma}[see \cite{Fretsaw} Lemma 8.14]\label{lem:freteigs}
Let $A,B$ be the stiffness matrices of $\mathcal{T},\mathcal{T}'$
respectively.  Let $B_S$ be the Schur complement of $B$ with respect to $A$.

If $\mathcal{T}'$ is a fretsaw extension of $\mathcal{T}$, 
then $\lambda_{min}(A,B_S) \geq 1/2$
\end{lemma}

We give a short proof of this lemma in Appendix \ref{apx:lemmas}.
Combining this with Lemma \ref{lambdamax}, we find:

\begin{corollary}\label{fretschur}
Let $A,B$ be the stiffness matrices of $\mathcal{T},\mathcal{T}'$
respectively.  Let $B_S$ be the Schur complement of $B$ with respect to $A$,
and let $A'=\lmx A & 0 \\ 0 & 0 \rmx$ be the same size as $B$.

If $\mathcal{T}'$ is a fretsaw extension of $\mathcal{T}$, 
then $\kappa(A,B_S) \leq 2\lambda_{max}(A',B)$
\end{corollary}

Now, consider a truss $\mathcal{T}$ and fretsaw extension $\mathcal{T}'$,
with respective rigidity graphs $Q_{\mathcal{T}}$ and $Q_{\mathcal{T}'}$.
By construction, every pair of faces that share an edge in $\mathcal{T}'$
must also share an edge in $\mathcal{T}$.  That is, if we let 
$\rho (Q_{\mathcal{T}'}) = \{(\rho (f_{1}),\rho (f_{2})):(f_{1},f_{2})\in Q_{\mathcal{T}'} \}$
denote the graph isomorphic to $Q_{\mathcal{T}'}$ on the faces of $\mathcal{T}$,
then $\rho (Q_{\mathcal{T}'})\subseteq Q_{\mathcal{T}}$.

As it turns out, for any subgraph $H\subseteq Q_{\mathcal{T}}$ of our choice,
we can construct a fretsaw extension with $Q_{\mathcal{T}'}$ (almost) 
isomorphic to $H$.
We present a linear-time construction here.  For technical reasons,
this construction also takes as input a map $\tau:[n]\rightarrow F$ that
for each vertex in $\mathcal{T}$ specifies one face containing that vertex.
The construction ensures that the face in $\mathcal{T}'$ corresponding to $\tau(i)$
contains the original copy of vertex $i$.
This feature will be useful later, and does not diminish
the generality of the algorithm.

\begin{lemma}
There exists an linear-time algorithm 
$\la \mathcal{T}',\rho \ra=\mathtt{fretsaw}(\mathcal{T},H,\tau)$
that takes a stiffly-connected truss 
$\mathcal{T}=\la n,\{\vv_{i} \}_{i=1}^{n},E,\gamma\ra$
with face set $F$,
a connected spanning 
subgraph $H$ of $Q_{\mathcal{T}}$, and a map $\tau:[n]\rightarrow F$ from
each vertex to a truss face containing it,
and returns a stiffly-connected 
$(\rho ,\pi )$-fretsaw extension 
$\mathcal{T}'=\la m,\{\vv_{i} \}_{i=1}^{m},E',\gamma'\ra \ra$
satisfying:
\begin{enumerate}
\item for all $i\in [n]$, $i \in \rho^{-1}(\tau(i))$
\item $H'\subseteq Q'$
\item
if $|H'|=n-1+k$ (i.e. $H'$ is a spanning tree plus $k$ additional edges),
then $|Q' - H'| \leq k$
\end{enumerate}
where $Q'$ denotes $Q_{\mathcal{T}'}$.
and $H'$ denotes $\rho^{-1}(H)=\{(f_{1},f_{2})\in Q':(\rho (f_{1}),\rho (f_{2}))\in H \}$
\end{lemma}

\begin{proof}
Here is the construction,
an example of which is given in Figure \ref{figure:truss}:

\mybox{
$\mathtt{fretsaw}(\mathcal{T},H,\tau)$

\medskip\noindent
First, for each vertex $i$ in $\mathcal{T}$,
we create the set $\pi^{-1}(i)$ of copies of vertex i
in $\mathcal{T}'$:

(Recall that we call $i\in \pi^{-1}(i)$ the ``original copy''.)
\begin{itemize}
\item 
Let $F_i$ denote the set of faces of $\mathcal{T}$ 
containing vertex $i$,
and let $H_i$ denote the graph induced by $H$ on $F_i$.
For each connected component of $H_{i}$, we put one copy 
of vertex $i$ in $\mathcal{T}'$.
The original copy is assigned to the connected component
of $H_{i}$ containing face $\tau (i)$.
\end{itemize}

Now, for $f\in F_{i}$, let $\phi (i,f)\in \pi^{-1}(i)$ 
denote the copy of vertex $i$
that is assigned to the component of $H_{i}$ containing $f$.
It is straighforward to construct the faces of $\mathcal{T}'$:
\begin{itemize}
\item For each face $f=(i,j,k)$ of $\mathcal{T}$,
we create a face $\rho^{-1} (f) = (\phi (i,f),\phi (j,f),\phi (k,f))$ in $\mathcal{T}'$.
\end{itemize}
}

The first property is directly enforced by the construction.

To see why the second property holds, consider an edge $(f_{1},f_{2})\in H$,
where $(i,j)$ is the edge shared by faces $f_{1}$ and $f_{2}$.  
Since edge $(f_{1},f_{2})$
is present in both $H_{i}$ and $H_{j}$, faces $\rho^{-1}(f_{1})$ and
$\rho^{-1}(f_{2})$ will share the same copies of vertices $i$ and $j$, and so
they too will share an edge.

As for the third property, suppose that $H'$ has $n-1+k$ edges, and thus
divides the plane into $k$ regions.
It suffices to show that each such region contains at most
one edge in $Q'-H'$.

Let $(f_{1},f_{2})$ be an edge in $Q_{\mathcal{T}'})-H'$.  
Let $(i,j)$ be the edge shared by $f_1$ and $f_2$.
Let $Q'_{i}$ and $H'_{i}\subseteq Q'_{i}$ denote the graphs induced on $F_{i}$ by 
$Q'$ and $H'$ respectively.

Since $f_{i}$ and $f_{j}$ share the vertices $i$ and $j$, 
we know there must a path from $f_{1}$ to $f_{2}$ both
in $H'_{i}$ and in $H'_{j}$.  Of course neither path contains
the edge $(f_{i},f_{2})$, since it is not in $H'$.
The only possibility then is that 
$Q'_{i}$ is a cycle $H'_{i}\cup \{(f_{1},f_{2}) \}$,
where $H'_{i}$ is a path from $f_{1}$ to $f_{2}$,
and similarly for $Q'_{j}$.
Thus $Q'_{i}\cup Q'_{j} = H'_{i}\cup H'_{j} \cup \{(f_{1},f_{2}) \}$,
and so $(f_{1},f_{2})$ is the only edge of $Q'$ inside the region
enclosed by cycle $H'_{i}\cup H'_{j}$.
\end{proof}

\section{Path Lemma}\label{pathsection}

We will need to construct
a fretsaw extension with a truss matrix that can
be solved quickly.  
In particular, the fretsaw extension we construct
will have a connectivity graph that is close to a spanning tree
(i.e. close to having $n-1$ edges),
because we can efficiently find a sparse 
Cholesky factorization of its truss matrix.
The following result is proven in Appendix (?somewhere?):

\begin{lemma}\label{cholesky_lemma}
Let $A$ be the stiffness matrix of an $n$-vertex truss $\mathcal{T}$, where
$Q_{\mathcal{T}}$ comprises a
spanning tree $R$ plus a set $S$ of 
additional edges.  A Cholesky factorization 
$A=PLL^TP^{T}$ can be found in time $O(n+|S|^{3/2})$, 
where $P$ is a permutation matrix, 
and such that systems in lower triangular matrix
$L$ can be used to solve systems in $A$ in time $O(n+|S|\log|S|)$.
\end{lemma}

Now, of course we want to construct a fretsaw extension 
whose truss matrix provides
good support for the original truss matrix.  
If we can give a supporting subset of faces
in the fretsaw extension for each
element in the original truss, then we can use Lemma \ref{cong-dil} to
bound the maximum generalized eigenvalue.
To this end, we show that a simply-connected set of faces that
connects two vertices supports the matrix of an element 
between the pair of vertices 
proportionally to the cube of the number of faces:

\begin{lemma}[Path Lemma]\label{lem:path}
Let $\mathcal{T}=\la n,\{\vv_{i} \}_{i=1}^{n},E,\gamma\ra$ be a $k$-face
simply-connected truss,
and let $e_{0}$ be a truss element with weight $\gamma(e_{0})$ between any pair
$\vv_{p},\vv_{q}\in V$ of vertices
in the truss, such that in $\mathcal{T}$ and $e_{0}$
\begin{itemize}
\item all truss elements have lengths in the range $[l_{min},l_{max}]$
\item all angles of truss faces are in the range $[\theta_{min},\pi-\theta_{min}]$.
\item all weights are in the range $[\gamma_{min},\gamma_{max}]$.
\end{itemize}
for positive constants $l_{min},l_{max},\theta_{min},\gamma_{min},\gamma_{max}$.
Then:
$$
\lambda_{max}(A_{e_{0}},A_{\mathcal{T}}) = O(k^3)
$$
\end{lemma}

We first note the simply-connected truss $\mathcal{T}$ must contain a 
simply-connected subset of faces whose rigidity graph is a path,
such that the first face in the path is the only one containing
vertex $p$ and the last triangle is the only one containing vertex $q$.
Thus, without loss of generality we may assume that $\mathcal{T}$ is
itself such a ``truss path'', because removing the extra faces can
only increase the value of $\lambda_{max}$.

Let us then number the faces in the path in order $f_1,...,f_n$,
and let us number the $n+2$ vertices as follows:
\begin{itemize}
\item 
Let $p=0$ be the vertex in $f_1$ but not $f_{2}$.
\item
Let vertices 1 and 2 be the pair of vertices shared by $f_{1}$ and $f_{2}$.
\item
For $3\leq i\leq n+1$, let $i$ be the vertex in $f_{i-1}$ but not $f_{i-2}$.
(In particular, $q=n+1$.)
\end{itemize}

\begin{figure}
\begin{center}
\epsfig{file=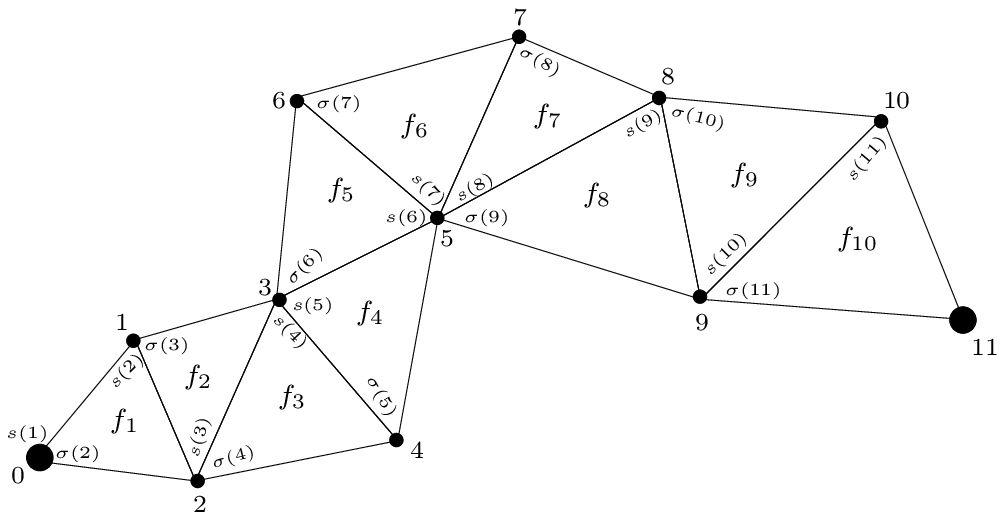, height=3in}
\end{center}
\caption{A truss path from $\vv_0$ to $\vv_{11}$, with triangles and
  vertices labeled appropriately.}
\label{trussfig}
\end{figure}
Furthermore, for $2\leq i\leq n$, consider the labeling of the three vertices in $f_{i-1}$
\begin{itemize}
\item $i$ is one of the vertices shared by $f_{i-1}$ and $f_{i}$.
\item
Let $s(i)$ denote the other vertex shared by $f_{i-1}$ and $f_{i}$.
\item
Let
$\sigma(i)$ denote the vertex in $f_{i-1}$ but not $f_i$.
\end{itemize}
For completeness, define $s(1)=0,s(n+1)=n,\sigma(n+1)=n-1$,
so for $1\leq i\leq n$:
\[
f_i = \{s(i),i,i+1\} = \{\sigma (i+1),s(i+1),i+1\}
\]
and the set of truss elements is given by:
\[
E = \{(0,1)\}\cup \bigcup_{i=2}^{n+1}\{(s(i),i),(\sigma(i),i) \}
\]
An example of this labeling is given in Figure \ref{trussfig}.

In the proof we make use of 
the following canonical definition of a perpendicular vector:
\begin{definition}
The {\bf counterclockwise perpendicular} of a vector 
$\xx = \lmx x_{1} \\ x_{2} \rmx \in \Reals2$ is
\[
 \xx^\perp = \lmx -x_{2} \\ x_{1} \rmx 
\]
\end{definition}
We note some useful properties of the perpendicular:
\begin{claim}\label{perpclaim}
For any $\xx,\yy,\zz\in\Reals2$:
\begin{enumerate}
\item $(-\xx)^\perp = -(\xx^\perp)$
\item $(\xx+\yy)^\perp = \xx^\perp + \yy^\perp$
\item $|\xx^T\yy^\perp|=|\xx||\yy||\sin\theta|$ where
  $\theta$ is the difference in angle between $\xx$ and $\yy$
\item
$
\frac{(\xx-\yy)^T(\xx-\zz)^\perp}
{(\xx-\yy)^T(\yy-\zz)^\perp} = 1
$
\item
$
\frac{\xx^\perp\yy^T}{\yy^T\xx^\perp} +
\frac{\yy^\perp\xx^T}{\xx^T\yy^\perp} = 
\lmx 1 & 0 \\ 0 & 1 \rmx = I
$
\end{enumerate}
\end{claim}

\begin{proof}
The first three properties are trivial.

The fourth can easily be seen by
$
\frac{(\xx-\yy)^T(\xx-\zz)^\perp}
{(\xx-\yy)^T(\yy-\zz)^\perp} = 
\frac{(\xx-\yy)^T(\xx-\yy)^\perp}
{(\xx-\yy)^T(\yy-\zz)^\perp} + 
\frac{(\xx-\yy)^T(\yy-\zz)^\perp}
{(\xx-\yy)^T(\yy-\zz)^\perp} =
0 + 1
$

Here is a proof of the fifth:
\[
\frac{\xx^\perp\yy^T}{\yy^T\xx^\perp} +
\frac{\yy^\perp\xx^T}{\xx^T\yy^\perp} = 
\frac{1}{x_{2}y_{1}-x_{1}y_{2}}
\lmx x_{2}y_{1} & x_{2}y_{2} \\
-x_{1}y_{1} & -x_{1}y_{2} \rmx -
\frac{1}{x_{2}y_{1}-x_{1}y_{2}}
\lmx x_{1}y_{2} & x_{2}y_{2} \\
-x_{1}y_{1} & -x_{2}y_{1} \rmx =
\lmx 1 & 0 \\ 0 & 1 \rmx = I
\]
\end{proof}

Now, let $\xx^*=[x^{*}_{0},x^{*}_{1},\dotsc,x^{*}_{2n-1}]^{T}$ be a vector that maximizes
$\frac{\xx^TA_{e_{0}}\xx}{\xx^TA_{\mathcal{T}}\xx}$
over all values of $\xx$.
Let $\xx^*_{i}$ denote $[x^{*}_{2i-1},x^{*}_{2i}]^{T}$

In particular we choose an $\xx^*$ such
that $(\xx^*_1-\xx^*_0)$ is parallel to $(\vv_1-\vv_0)$,
by taking advantage of the following property of the null space:
\begin{lemma}
Define $\xx^R=[x^{R}_{0},\dotsc ,x^{R}_{2n-1}]^{T}$ to be the vector satisfying
$[x^{R}_{2i-1},x^{R}_{2i}]^{T} = (\vv_i-\vv_0)^\perp $. $\xx^R$ is in the null space of both
$A_{e_{0}}$ and $A_{\mathcal{T}}$.
\end{lemma}
\begin{proof}
For the matrix $A_{e}$ of any single element $e=(i,j)$, we have:
\[
A_{e}\xx^{R} =
\frac{\gamma}{|\vv_i-\vv_j|} \uu_e \uu_e^T\xx^{R} =
\frac{\gamma}{|\vv_i-\vv_j|^{3}} \uu_e 
\left((\vv_{i}-\vv_{j})^T(\vv_i-\vv_j)^\perp + 
(\vv_{j}-\vv_{i})^T(\vv_j-\vv_i)^\perp\right) = 0
\]
\end{proof}
Note that we can eliminate the component of $\xx^{*}$ perpendicular to
$(\vv_{1}-\vv_{0})$ by adding the appropriate multiple of $\xx^{R}$.

Now, let us focus momentarily on a single vertex $i$, and the two
elements $(s(i),i)$ and $(\sigma (i),i)$ that connect vertex $i$ to
lower numbered vertices.  The terms $\xx^{T}A_{(s(i),i)}\xx$ and
$\xx^{T}A_{(\sigma (i),i)}\xx$ are zero respectively when 
\[
(\vv_{i}-\vv_{s(i)})^T(\xx_{i}-\xx_{s(i)})=0
\]
and
\[
(\vv_{i}-\vv_{\sigma (i)})^T(\xx_{i}-\xx_{\sigma (i)})=0
\]
Supposing we set $\xx_{s(i)}=\xx^{*}_{s(i)}$ and
$\xx_{\sigma (i)}=\xx^{*}_{\sigma (i)}$, we would like to define
$\dd_{i}$ to be the vector such that setting $\xx_{i}=\xx^{*}_{i}-\dd_{i}$
satisfies both of the above equations.

In particular, we define the vectors
$$
\dd_i = 
\begin{cases}
\xx^*_1 - \xx^*_0 & i = 1 \\
\xx^*_i - \xx^*_{s(i)} + R_i(\xx^*_{s(i)}-\xx^*_{\sigma(i)})
& 2 \leq i \leq n+1 \\
\end{cases}
$$
where
$$
R_i = 
\frac{(\vv_i - \vv_{s(i)})^\perp(\vv_i-\vv_{\sigma(i)})^T}
{(\vv_i-\vv_{\sigma(i)})^T(\vv_i-\vv_{s(i)})^\perp}
$$

We claim that these satisfy the following properties:
\begin{lemma}\label{lem:dprops}
The following are properties of the $\dd_i$s:
\begin{enumerate}
\item 
For all $(j,i)\in E$, $j<i$:
\[
(\vv_i-\vv_{j})^T(\xx^*_i-\xx^*_{j}) = 
(\vv_i-\vv_{j})^T\dd_i
\]
\item
For all $i$:
\[
|\xx^*_i-\xx^*_{s(i)}| \leq
\frac{l_{max}}{l_{min}\sin\theta_{min}}\sum_{j=1}^i|\dd_j|
\]
\end{enumerate}
\end{lemma}

\begin{proof}[Proof of 1]
The statement is trivial for element $(0,1)$
There are two other types of elements we must consider: $(s(i),i)$ and $(\sigma (i),i)$.

For an element $(s(i),i)$ we have:
\begin{align*}
(\vv_{i}-\vv_{s(i)})^{T}\dd_{i} 
&=
(\vv_{i}-\vv_{s(i)})^{T}(\xx^{*}_{i}-\xx^{*}_{s(i)}) +
\frac{(\vv_{i}-\vv_{s(i)})^{T}(\vv_i - \vv_{s(i)})^\perp
(\vv_i-\vv_{\sigma(i)})^T(\xx^*_{s(i)}-\xx^*_{\sigma(i)})}
{(\vv_i-\vv_{\sigma(i)})^T(\vv_i-\vv_{s(i)})^\perp} \\
&= (\vv_{i}-\vv_{s(i)})^{T}(\xx^{*}_{i}-\xx^{*}_{s(i)}) + 0
\intertext{using the fact that $(\vv_{i}-\vv_{s(i)})^{T}(\vv_i - \vv_{s(i)})^\perp = 0$}
\end{align*}

For an element $(\sigma (i),i)$ we have:
\begin{align*}
(\vv_{i}-\vv_{\sigma (i)})^{T}\dd_{i} 
&=
(\vv_{i}-\vv_{\sigma (i)})^{T}(\xx^{*}_{i}-\xx^{*}_{s(i)}) +
\frac{(\vv_{i}-\vv_{\sigma (i)})^{T}(\vv_i - \vv_{s(i)})^\perp
(\vv_i-\vv_{\sigma(i)})^T(\xx^*_{s(i)}-\xx^*_{\sigma(i)})}
{(\vv_i-\vv_{\sigma(i)})^T(\vv_i-\vv_{s(i)})^\perp} \\
&= (\vv_{i}-\vv_{\sigma (i)})^{T}(\xx^{*}_{i}-\xx^{*}_{s(i)}) +
(\vv_i-\vv_{\sigma(i)})^T(\xx^*_{s(i)}-\xx^*_{\sigma(i)}) \\
&= (\vv_{i}-\vv_{\sigma (i)})^{T}(\xx^{*}_{i}-\xx^{*}_{\sigma (i)})
\end{align*}
\end{proof}
\begin{proof}[Proof of 2]
For $i\geq 2$, using the fact that $\{s(i),\sigma (i) \}=\{s(i-1),i-1 \}$,
we have
$$
\xx^*_i-\xx^*_{s(i)} =
\dd_i - R_i(\xx^*_{s(i)}-\vv^*_{\sigma(i)}) =
\dd_i \pm R_i(\xx^*_{i-1}-\vv^*_{s(i-1)})
$$

Since $\xx^*_1 - \xx^*_0 = \dd_1$, we recursively find that
$$
|\xx^*_i-\xx^*_{s(i)}| \leq 
|\dd_i| + \sum_{j=1}^{i-1}
|(R_iR_{i-1}\cdots R_{j+1})\dd_j|
$$

The following finishes the proof:
\begin{align*}
& |R_iR_{i-1}\cdots R_{j+1}\dd_j| \\
&=
\left|
(\vv_i-\vv_{s(i)})^\perp
\prod_{k=j+2}^i \left(
\frac{(\vv_k-\vv_{\sigma(k)})^T(\vv_{k-1}-\vv_{s(k-1)})^\perp}
{(\vv_k-\vv_{\sigma(k)})^T(\vv_k-\vv_{s(k)})^\perp} \right)
\cdot
\frac{(\vv_{j+1}-\vv_{\sigma(j+1)})^T\dd_j}
{(\vv_{j+1}-\vv_{\sigma(j+1)})^T(\vv_{j+1}-\vv_{s(j+1)})^\perp} 
\right| \\
&=
\left|
(\vv_i-\vv_{s(i)})^\perp
\prod_{k=j+2}^i \left(
\frac{(\vv_k-\vv_{\sigma(k)})^T(\vv_{s(k)}-\vv_{\sigma(k)})^\perp}
{(\vv_k-\vv_{\sigma(k)})^T(\vv_k-\vv_{s(k)})^\perp} \right)
\cdot
\frac{(\vv_{j+1}-\vv_{\sigma(j+1)})^T\dd_j}
{(\vv_{j+1}-\vv_{\sigma(j+1)})^T(\vv_{j+1}-\vv_{s(j+1)})^\perp} 
\right| \\
\intertext{by the fact that $\{s(i),\sigma (i) \}=\{s(i-1),i-1 \}$}
&=
\left|
(\vv_i-\vv_{s(i)})^\perp
\prod_{k=j+2}^i \left( \pm 1 \right)
\cdot
\frac{(\vv_{j+1}-\vv_{\sigma(j+1)})^T\dd_j}
{(\vv_{j+1}-\vv_{\sigma(j+1)})^T(\vv_{j+1}-\vv_{s(j+1)})^\perp} 
\right| \\
\intertext{by Claim \ref{perpclaim}.4}
&=
\left|
\frac{(\vv_i-\vv_{s(i)})^\perp(\vv_{j+1}-\vv_{\sigma(j+1)})^T
\dd_j}
{(\vv_{j+1}-\vv_{\sigma(j+1)})^T(\vv_{j+1}-\vv_{s(j+1)})^\perp} 
\right| \\
&\leq 
\frac{|\vv_i-\vv_{s(i)}||\vv_{j+1}-\vv_{\sigma(j+1)}||\dd_j|}
{|\vv_{j+1}-\vv_{\sigma(j+1)}||\vv_{j+1}-\vv_{s(j+1)}|\sin \theta_{min}} 
\\
\intertext{by Claim \ref{perpclaim}.3}
&=
\frac{|\vv_i-\vv_{s(i)}||\dd_j|}
{|\vv_{j+1}-\vv_{s(j+1)}|\sin \theta_{min}}\\
&\leq
\frac{l_{max}}{l_{min}\sin\theta_{min}}|\dd_j| \\
\end{align*}
\end{proof}

To finish proving the path lemma, we will need to
use the following fact:
\begin{lemma}\label{lem:2vecs}
Let $\uu_1,\uu_2$ be unit vectors whose angles differ by $\theta$.
Then for any $\vv$ and $a,b>0$:
$$
\frac{(\uu_1^T\vv)^2}{a} + \frac{(\uu_2^T\vv)^2}{b} \geq
\frac{\sin^2\theta|\vv|^2}{a+b}
$$
\end{lemma}
\begin{proof}
Let $\alpha$ be the angle between $\uu_{1}$ and $\vv$.

We must show that 
$\frac{1}{a}\cos^{2}\alpha + \frac{1}{b}\cos^{2}(\alpha +\theta ) \geq 
\frac{1}{a+b}\sin^{2}\theta$
\end{proof}

Recall that we wish to prove 
\[
\lambda_{max}(A_{e_{0}},A_{\mathcal{T}})
=
\frac{(\xx^*)^TA_{e_{0}}\xx^*}{(\xx^*)^TA_{\mathcal{T}}\xx^*}
= O(n^{3})
\]

Let us first bound the denominator $(\xx^*)^TA_{\mathcal{T}}\xx^*$.
\begin{align*}
& (\xx^*)^TA_{\mathcal{T}}\xx^* \\
&=
(\xx^{*})^TA_{(0,1)}\xx^{*} +
\sum_{i=2}^{n+1}
\left((\xx^{*})^TA_{(s(i),i)}\xx^{*} +
(\xx^{*})^TA_{(\sigma (i),i)}\xx^{*} \right) \\
&\geq 
\gamma_{min}\left[
\frac{\left[\left(\frac{\vv_1-\vv_0}{|\vv_1-\vv_0|}\right)^T
(\xx^*_1-\xx^*_0)\right]^2}
{|\vv_1-\vv_0|} + 
\sum_{i=2}^{n+1}
\left( 
\frac{\left[\left(\frac{\vv_i-\vv_{s(i)}}
{|\vv_i-\vv_{s(i)}|}\right)^T
(\xx^*_i-\xx^*_{s(i)})\right]^2}
{|\vv_i-\vv_{s(i)}|} + 
\frac{\left[\left(\frac{\vv_i-\vv_{\sigma(i)}}
{|\vv_i-\vv_{\sigma(i)}|}\right)^T
(\xx^*_i-\xx^*_{\sigma(i)})\right]^2}
{|\vv_i-\vv_{\sigma(i)}|} \right) 
\right] \\
&=
\gamma_{min}\left[
\frac{|\dd_{1}|^2}
{|\vv_1-\vv_0|} + 
\sum_{i=2}^{n+1}
\left( 
\frac{\left[\left(\frac{\vv_i-\vv_{s(i)}}
{|\vv_i-\vv_{s(i)}|}\right)^T
\dd_i\right]^2}
{|\vv_i-\vv_{s(i)}|} + 
\frac{\left[\left(\frac{\vv_i-\vv_{\sigma(i)}}
{|\vv_i-\vv_{\sigma(i)}|}\right)^T
\dd_i\right]^2}
{|\vv_i-\vv_{\sigma(i)}|} \right) 
\right] \\
\intertext{by Lemma \ref{lem:dprops} and the fact that
$\dd_{1}=\xx^{*}_{1}-\xx^{*}_{0}$ is parallel to $(\vv_{1}-\vv_{0})$}
&\geq
\gamma_{min}\left[
\frac{|\dd_1|^2}
{|\vv_1-\vv_0|} + 
\sum_{i=2}^{n+1}
\left( 
\frac{\sin^2\theta_{min}|\dd_i|^2}
{|\vv_i-\vv_{s(i)}|+
|\vv_i-\vv_{\sigma(i)}|}
\right) 
\right] \\
\intertext{by Lemma \ref{lem:2vecs}}
&\geq
\sum_{i=1}^{n+1}
\left( 
\frac{\gamma_{min}\sin^2\theta_{min}|\dd_i|^2}
{2l_{max}}
\right) \\
\end{align*}

Next we bound the numerator:

\begin{align*}
(\xx^*)^TA_{e_{0}}\xx^*
&=
\frac{\gamma (e_{0})}{|\vv_{n+1}-\vv_{0}|}
\left[ \left( \frac{\vv_{n+1}-\vv_0}{|\vv_{n+1}-\vv_0|} \right)^T
(\xx^*_{n+1}-\xx^*_0) \right]^2 \\
&\leq
\frac{\gamma_{max}}{l_{min}}
|\xx^*_{n+1}-\xx^*_0|^2 \\
&\leq 
\frac{\gamma_{max}}{l_{min}}
\left(
|\xx^*_{n+1}-\xx^*_{s(n+1)}| +
|\xx^*_{s(n+1)}-\xx^*_{s(s(n+1))}| + \cdots +
|\xx^*_1-\xx^*_0|
\right)^2 \\
&\leq
\frac{\gamma_{max}}{l_{min}}
\left(
(n+1)\frac{l_{max}}{l_{min}\sin\theta_{min}}\sum_{i=1}^{n+1}|\dd_i|
\right)^2
\intertext{by Lemma \ref{lem:dprops}}
\end{align*}

Combining the above, we get:

\begin{align*}
\lambda_{max}(A_{e_{0}},A_{\mathcal{T}})
&=
\frac{(\xx^*)^TA_{e_{0}}\xx^*}{(\xx^*)^TA_{\mathcal{T}}\xx^*} \\
&\leq
\frac{2}{\sin^{3}\theta_{min}}
\left(\frac{l_{min}}{l_{max}} \right)^{3}
\frac{\gamma_{max}}{\gamma_{min}}
\frac{\left(\sum_{i=1}^{n+1}|\dd_i|\right)^2}
{\sum_{i=1}^{n+1}|\dd_i|^2} 
(n+1)^2 \\
&\leq
\frac{2}{\sin^{3}\theta_{min}}
\left(\frac{l_{min}}{l_{max}} \right)^{3}
\frac{\gamma_{max}}{\gamma_{min}}
(n+1)^3 \\
\end{align*}
where the last inequality follows from Cauchy-Schwarz.

\section{Graph Embeddings}\label{stretch}

Our remaining task is to describe how to map edges to
supporting face sets with low congestion, as required for
Lemma \ref{cong-dil}.
We need some graph theoretic notions which will inform how we choose
to support
the edges.  

For a graph with edges $E$, we use $\mathcal{P}(E)$ to denote the power
set of $E$.  Thus $\mathcal{P}(E)$ includes all paths in the graph.

Let us define the notion of embedding 
vertex pairs of a graph into paths in a subgraph:

\begin{definition}
For an unweighted graph $(V,E)$, a subgraph $H\subseteq E$, 
and a set $Z$ of pairs of vertices in $V$,
an {\bf embedding} of $Z$ onto $H$ is 
a map $\pi: Z \rightarrow \mathcal{P}(H)$, where
$\pi(v,w)\subseteq H$ is a path in $H$ whose endpoints are $v,w$.

The {\bf stretch} of $\pi$ is 
$str(\pi) = \sum_{z\in Z} |\pi(z)|$.

The {\bf congestion} of $\pi$ is 
$cong(\pi) = \max_{f\in H} \left[ \sum_{z\in Z:f\in \pi(z)} |\pi(z)| \right]$
\end{definition}
A particular example is embedding the edges of a graph onto a spanning tree:
\begin{definition}
For an unweighted connected graph $(V,E)$ and spanning tree $T\subseteq E$, 
let $T(v,w)$ denote the path in $T$ that connects $v$ to $w$.

The {\bf stretch} of $T$ is 
$str(T) = \sum_{e\in E} |T(e)|$.
\end{definition}

We will make use of algorithms that take a graph and generate a spanning tree,
augmented with a few additional edges, such that given vertex pairs have
a low congestion embedding into the augmented tree.
First, we need to create a low stretch-spanning tree.  The best known result
is from \cite{EEST}:

\begin{theorem}
There exists an algorithm $T=\mathtt{LowStretch}(G)$,
that takes a connected graph $G=(V,E)$, runs in time $O(|E|\log^2 |E|)$,
and outputs a spanning tree $T$ with stretch 
$O(|E|\log^2 |V|\log\log |V|)$.
\end{theorem}

We can then use the low-stretch spanning tree to create an augmented spanning
tree with the desired low congestion embedding.  This algorithm is given in 
Appendix \ref{congapp}:

\begin{theorem}\label{lowcong}
There exists an algorithm $S=\mathtt{LowCongestAugment}(G,T,Z,\psi,k)$
that takes a planar graph $G=(V,E)$, a spanning tree $T$ of $G$,
a set $Z$ of pairs of vertices in $V$, an embedding 
$\psi: Z\rightarrow\mathcal{P}(E)$, and an integer $k$.
The algorithm runs in time $O(|E|\log |V| + cong(\psi)|E|)$
and returns a set of edges $S\subseteq E$ of
size at most $k$,
such that there exists an embedding $\pi:Z\rightarrow\mathcal{P}(T\cup S)$
with congestion $O(\frac{1}{k}str(T)cong(\psi))$.
\end{theorem}

\section{Solving the Linear System}

In Figure \ref{mainalg} we present
the complete {\tt TrussSolver} 
algorithm for solving linear systems in matrix a $A$ that is
the stiffness matrix of truss 
$\mathcal{T} = \la n,\{\vv_{i} \}_{i=1}^{n},E,\gamma\ra$
with face set $F$.
The algorithm preconditions $A$ using the stiffness matrix $B$
of a fretsaw extension
$\mathcal{T}' = \la m,\{\vv'_{i} \}_{i=1}^{m},E',\gamma'\ra$
with face set $F'$.
The algorithm uses a parameter $k$ that will be chosen later.
We will show that (with the right choice of $k$) the algorithm attains
a relative error of $\epsilon$ in time 
$O\left(n^{5/4}(\log^2 n\log\log n)^{3/4}\log\frac{1}{\epsilon}\right)$.

\begin{figure}
\mybox{
$\mathtt{TrussSolver}(\mathcal{T},k)$

\medskip\noindent
Let $A$ be the stiffness matrix of $\mathcal{T}=\la n,\{\vv_{i} \}_{i=1}^{n},E,\gamma\ra$,
and let $F$ be its face set.
\begin{enumerate}
\item
Define $\tau:[n]\rightarrow F$ to map vertex 
$i$ to an arbitrary face in 
the set $F_i$ of faces containing $i$.

\item
Run $R = \mathtt{LowStretch}(Q_{\mathcal{T}})$ in time $O(n\log^2 n)$.

\item
Define $Z = \{(\tau(i),\tau(j)):(i,j)\in E\}$.

Define an embedding $\psi: Z\rightarrow\mathcal{P}(Q_{\mathcal{T}})$, by defining
$\psi(\tau(i),\tau(j))$ to be an arbitrary path from $\tau(i)$ to $\tau(j)$ in
$F_i \cup F_j$.  (We know $F_i$ and $F_j$ intersect because some face contains
edge $(i,j)$.)

Run $S = \mathtt{LowCongestAugment}(Q_{\mathcal{T}},R,Z,\psi,k)$ in time $O(n\log n)$.

\item
Run $\la \mathcal{T}',\rho \ra = 
\la \la V',E',\gamma' \ra,\rho \ra = 
\mathtt{fretsaw}(\mathcal{T},R\cup S,\tau)$.
Let $B$ be the truss matrix of $\mathcal{T}'$.

\item
Use Lemma \ref{cholesky_lemma} to
find a Cholesky factorization $B = PLL^TP^{T}$ in time $O(n+k^{3/2})$,
such that $L$ can be used to solve equations in $B$ in time $O(n+k\log k)$.

\item
Run preconditioned conjugate gradient using $B_S$, the Schur complement of $B$ with
respect to $A$, as the preconditioner. 
Use $L$ to solve equations in $B_S$, by solving equations in $B$ (see Lemma
\ref{equivlem}).

The relative error will be down to $\epsilon$ after
$O(\sqrt{\kappa(A,B_S)}\log\frac{1}{\epsilon})$ iterations.
\end{enumerate}
}
\caption{The {\tt TrussSolver} algorithm}
\label{mainalg}
\end{figure}

We want to use the Congestion Dilation Lemma (Lemma \ref{cong-dil}) 
to give an upper bound on $\kappa(A,B_S)$.
Recall that for $A' = \lmx A & 0 \\ 0 & 0 \rmx$ of the same size as $B$,
Corollary \ref{fretschur} says that $\kappa(A,B_S)\leq 2\lambda_{max}(A',B)$.
To bound $\lambda_{max}(A',B)$, 
for each edge $(p,q)\in E$ we will
give a face subset $F'_{p,q}\subseteq F'$ connecting $p$ to $q$,
such that this
embedding of edges to truss paths has low congestion.

In particular, let $E'_{p,q}\subseteq E'$ be the set of truss elements in
the faces $F_{p,q}$, so that
$\mathcal{T}'_{p,q}=\la m,\{\vv'_{i}\},E'_{p,q},\gamma'\ra$ 
is the ``subtruss'' of $\mathcal{T}'$
comprising the faces $F_{p,q}$.
Lemma \ref{lem:path} states that
$$ \lambda_{max}(A_{(p,q)},A_{\mathcal{T}_{p,q}}) = O(|F_{p,q}|^3) $$
so Lemma \ref{cong-dil} yields
$$
\lambda_{max}(A',B) 
\leq
\max_{e\in E'} \sum_{(p,q)\in E:e \in E'_{p,q}} O(|F'_{p,q}|^3)
$$
It remains for us to describe the truss paths $\mathcal{T}'_{p,q}$ 
that yield the desired bound.

Recall that we have constructed a subgraph $R\cup S\in Q_{\mathcal{T}}$,
for which the {\tt LowCongestAugment} algorithm guarantees that
there exists an embedding
$\pi:Z\rightarrow\mathcal{P}(R\cup S)$
of low congestion.
Let us denote $\pi_{p,q} = \pi(\tau(p),\tau(q))$, the path in 
$R\cup S\subseteq Q_{\mathcal{T}}$
connecting $\tau(p)$ to $\tau(q)$.  Map this path back into 
$Q_{\mathcal{T}'}$ to get the path
$\pi'_{p,q} = \rho^{-1}(\pi_{p,q}) = 
\{(\rho^{-1}(f_1),\rho^{-1}(f_2)):(f_{1},f_{2})\in\pi_{p,q}\}$.
We then form $\mathcal{T}'_{p,q}$ from the set
$F'_{p,q}$ of faces in $\pi'_{p,q}$.

Let us first determine the congestion of $\pi$ more precisely.
The algorithms $\mathtt{LowCongestAugment}$ and $\mathtt{LowStretch}$ 
guarantee respectively that
$cong(\pi) = O(\frac{1}{k}str(R)cong(\psi))$
and
$str(R) = O(n\log^2 n \log\log n)$.

As for $cong(\psi)$, we have 
$$
cong(\psi) = 
\max_{q\in Q_{\mathcal{T}}} \sum_{z\in Z:q\in\psi(z)}|\psi(z)| \leq
\left(\max_{z\in Z} |\psi(z)| \right)
\left(\max_{q\in Q} |\{z\in Z:q\in\psi(z)\}| \right)
$$
Now note that for all $i$, 
$|F_i|\leq \frac{2\pi}{\theta_{min}}=O(1)$.  
Since any $\psi(\tau_i,\tau_j)$ only contains triangles in $T_i\cup T_j$,
$|\psi(\tau(i),\tau(j))|=O(1)$, and so
$\max_{z\in Z} |\psi(z)|=O(1)$.

Similarly, say that $q\in Q_{\mathcal{T}}$ 
is a pair of faces sharing the edge
$(i,j)$.  Since $i$ and $j$ are the only vertices that the pair of faces
have in common, $q$ can only be in a path $\psi(\tau_\alpha, \tau_\beta)$ if
one of $\alpha$ or $\beta$ is $i$ or $j$.
So $|\{z:q\in\psi(z)\}|\leq
|\{(\alpha,\beta)\in E:\alpha=i\text{ or }\alpha=j\}|\leq
\frac{4\pi}{\theta_{min}}=O(1)$.

Thus, $cong(\psi) = O(1)$ and $cong(\pi)=O(\frac{n}{k}\log^2 n\log\log n)$.

We now have:
\begin{eqnarray*}
\kappa(A,B_S) 
&\leq&
\max_{e\in E'} \sum_{(p,q)\in E:e \in E'_{p,q}} O(|F'_{p,q}|^3) \\
&=&
O\left(\max_{e\in E'} \mathop{\sum_{(p,q)\in E:}}_{e \in E'_{p,q}}
|F'_{p,q}|\right)^3 \\
&=&
O\left(\max_{(f_i,f_j)\in Q_{\mathcal{T}'}} 
\mathop{\sum_{(p,q)\in E:}}_{(f_i,f_j) \in \pi'_{p,q}}
|\pi'_{p,q}| \right)^3 \\
&=&
O\left(\max_{(f_i,f_j)\in R\cup S} 
\mathop{\sum_{(p,q)\in E:}}_{(f_i,f_j) \in \pi_{p,q}}
|\pi_{p,q}| \right)^3 \\
&=&
O\left(\max_{(f_i,f_j)\in R\cup S} 
\mathop{\sum_{z\in Z:}}_{(f_i,f_j) \in \pi(z)}
|\pi(z)| \right)^3 \\
&=&
O\left(cong(\pi)\right)^3 \\
&=&
O\left( \frac{n^3}{k^3}(\log^2 n \log\log n)^3 \right) \\
\end{eqnarray*}

Steps 1-4 take time $O(n\log^2 n + k^{3/2})$.
Each conjugate gradient iteration takes time $O(n+k\log k)$, and the number of
iterations is
$$ 
\sqrt{\kappa(A,B_S)}\log\frac{1}{\epsilon} =
O(\frac{n^{3/2}}{k^{3/2}}(\log^2 n\log\log n)^{3/2})
\log\frac{1}{\epsilon}
$$
Thus, our total running time is:
$$
O\left(n\log^2 n+k^{3/2} +n^{3/2}k^{-3/2}(\log^2 n\log\log n)^{3/2}(n+k\log k)
\log\frac{1}{\epsilon}\right)
$$
For $k=n^{5/6}(\log^2 n\log\log n)^{1/2}$ 
this gives a running time of 
$O\left(n^{5/4}(\log^2 n\log\log n)^{3/4}\log\frac{1}{\epsilon}\right)$.

\bibliography{precon}

\begin{thebibliography}{EEST06}

\bibitem[Axe85]{Axelsson}
O.~Axelsson.
\newblock A survey of preconditioned iterative methods for linear systems of
  algebraic equations.
\newblock {\em BIT Numerical Mathematics}, 25(1):165--187, March 1985.

\bibitem[BH03]{SupportTheory}
Erik~G. Boman and Bruce Hendrickson.
\newblock Support theory for preconditioning.
\newblock {\em SIAM Journal on Matrix Analysis and Applications},
  25(3):694--717, 2003.

\bibitem[EEST06]{EEST}
Michael Elkin, Yuval Emek, Daniel~A. Spielman, and Shang-Hua Teng.
\newblock Lower-stretch spanning trees.
\newblock {\em SIAM Journal on Computing}, 2006.
\newblock To appear.

\bibitem[ST06a]{Fretsaw}
Gil Shklarski and Sivan Toledo.
\newblock Rigidity in finite-element matrices: Sufficient conditions for the
  rigidity of structures and substructures, 2006.

\bibitem[ST06b]{SpielmanTengLinsolve}
Daniel~A. Spielman and Shang-Hua Teng.
\newblock Nearly-linear time algorithms for preconditioning and solving
  symmetric, diagonally dominant linear systems.
\newblock Available at \texttt{http://www.arxiv.org/abs/cs.NA/0607105}, 2006.

\end{thebibliography}
\bibliographystyle{alpha}

\appendix

\section{Preconditioning Lemmas}\label{apx:lemmas}

We first prove several lemmas dealing with the Schur complement.
Recall that the Schur complement of 
$B = \lmx B_{11} & B_{12} \\ B_{12}^T & B_{22} \rmx$
is
$B_{S} = B_{11} - B_{12}B_{22}^{-1}B_{12}^T$

\begin{lemma}\label{schurprop}
If $B$ is positive semidefinite then for any $\xx$
$$
\min_{\yy} \left(
\lmx \xx^T & \yy^T \rmx B \lmx \xx \\ \yy \rmx \right) =
\xx^T B_S \xx
$$
\end{lemma}

\begin{proof}
\begin{eqnarray*}
\lmx \xx^T & \yy^T \rmx B \lmx \xx \\ \yy \rmx
&=&
\xx^TB_{11}\xx +
\xx^TB_{12}\yy +
\yy^TB_{12}^T\xx +
\yy^TB_{22}\yy \\
&=&
\xx^T(B_{11}-B_{12}B_{22}^{-1}B_{12}^T)\xx +
\xx^TB_{12}B_{22}^{-1}B_{12}^T\xx +
\xx^TB_{12}\yy +
\yy^TB_{12}^T\xx +
\yy^TB_{22}\yy \\
&=&
\xx^TB_S\xx +
(\yy+B_{22}^{-1}B_{12}^T\xx)^TB_{22}(\yy+B_{22}^{-1}B_{12}^T\xx) \\
&=&
\xx^TB_S\xx +
\lmx 0 \\ \yy+B_{22}^{-1}B_{12}^T\xx \rmx^T
B
\lmx 0 \\ \yy+B_{22}^{-1}B_{12}^T\xx \rmx \\
&\geq&
\xx^TB_S\xx \\
\end{eqnarray*}
The last inequality holds because $B$ is positive semidefinite,
and it is an equality for
$\yy = -B_{22}^{-1}B_{12}^T\xx$.
\end{proof}

\newtheorem*{lemmaC}{Lemma \ref{lambdamax}}
\begin{lemmaC}
For positive semidefinite $A,B$,
$\lambda_{max}(A',B) = \lambda_{max}(A,B_S)$.
\end{lemmaC}
\begin{proof}
\begin{eqnarray*}
\lambda_{max}(A',B)
&=&
\max_{\xx,\yy}
\frac{\lmx \xx^T & \yy^T \rmx A' \lmx \xx \\ \yy \rmx}
{\lmx \xx^T & \yy^T \rmx B \lmx \xx \\ \yy \rmx} \\
&=&
\max_{\xx,\yy}
\frac{\xx^T A \xx}
{\lmx \xx^T & \yy^T \rmx B \lmx \xx \\ \yy \rmx} \\
&=&
\max_{\xx}
\frac{\xx^TA\xx}{\xx^TB_S\xx} \\
&=&
\lambda_{max}(A,B_S) \\
\end{eqnarray*}
The third equality uses Lemma \ref{schurprop}.
\end{proof}

\newtheorem*{lemmaA}{Lemma \ref{equivlem}}
\begin{lemmaA}
$B\lmx \xx \\ \yy \rmx = 
\lmx \bb \\ 0 \rmx$ implies
$B_S\xx = \bb$
\end{lemmaA}

\begin{proof}
Multiplying by 
$\lmx I & -B_{12}B_{22}^{-1} \\ 0 & I \rmx$
gives:
\begin{eqnarray*}
\lmx I & -B_{12}B_{22}^{-1} \\ 0 & I \rmx
\lmx B_{11} & B_{12} \\ B_{12}^T & B_{22} \rmx 
\lmx \xx \\ \yy \rmx 
&=& 
\lmx I & -B_{12}B_{22}^{-1} \\ 0 & I \rmx
\lmx \bb \\ 0 \rmx \\
\lmx B_S & 0 \\ B_{12}^T & B_{22} \rmx 
\lmx \xx \\ \yy \rmx 
&=&
\lmx \bb \\ 0 \rmx \\
\end{eqnarray*}
\end{proof}

\newtheorem*{lemmaD}{Lemma \ref{cong-dil}}
\begin{lemmaD}[Congestion-Dilation Lemma]

Given the symmetric positive semidefinite matrices
$A_1,...,A_n,B_1,...,B_m$ and $A = \sum_i A_i$ and $B = \sum_j B_j$
and given sets $\Sigma_i \subseteq [1,...,m]$ and real values $s_i$
that satisfy

$$
\lambda_{max}(A_i,\sum_{j\in \Sigma_i} B_j) \leq s_i
$$

it holds that

$$
\lambda_{max}(A,B) \leq \max_j \left(\sum_{i:j\in\Sigma_i} s_i\right)
$$

\end{lemmaD}

\comment{
We will need to use the following lemma (see \cite{Gremban} Lemma 4.7):

\begin{lemma}\label{splitting}
Given the symmetric positive semidefinite matrices
$A_1,...,A_n,B_1,...,B_n$,
$$
\lambda_{max}(\sum_i A_i, \sum_i B_i) \leq \max_i (A_i,B_i)
$$
\end{lemma}
}

\begin{proof}
Let us write $A \cleq B$ to mean that $B-A$ is positive semidefinite.
We are given that 
\[
A_{i} \cleq s_{i}\sum_{j\in \Sigma_{i}}B_{j}
\]
So we have:
\begin{align*}
A &= \sum_{i}A_{i} \\
&\cleq \sum_{i}s_{i}\sum_{j\in \Sigma_{i}}B_{j}\\
&= \sum_{j}B_{j}\sum_{i:j\in \Sigma_{i}}s_{i}\\
&\cleq \max_{j}\left(\sum_{i:j\in \Sigma_{i}}s_{i}\right)\sum_{j}B_{j}\\
&= \max_{j}\left(\sum_{i:j\in \Sigma_{i}}s_{i}\right)B\\
\end{align*}

\comment{
\begin{eqnarray*}
\lambda_{max}(A,B)
&=&
c_{max} \lambda_{max}(A,c_{max} B) \\
(*)
&\leq&
c_{max} \lambda_{max}(A,\sum_j c_j B_j) \\
&=&
c_{max} \lambda_{max}(A,\sum_j \sum_{i:j\in\Sigma_i} s_i B_j) \\
&=&
c_{max} \lambda_{max}(\sum_i A_i,\sum_i s_i \sum_{j\in\Sigma_i} B_j) \\
(**)
&\leq&
c_{max} \max_i \lambda_{max}(A_i, s_i \sum_{j\in\Sigma_i} B_j) \\
&=&
c_{max} \max_i \frac{1}{s_i}\lambda_{max}(A_i, \sum_{j\in\Sigma_i} B_j) \\
&=&
c_{max} \\
\end{eqnarray*}
Inequality (*) follows from the fact that
$c_{max}B - \sum_j c_j B_j = \sum_j (c_{max}-c_j)B_j$
is positive semidefinite
and inequality (**) uses Lemma \ref{splitting}.
}
\end{proof}

\newtheorem*{lemmaE}{Lemma \ref{lem:freteigs}}
\begin{lemmaE}[see \cite{Fretsaw} Lemma 8.14]
Let $A,B$ be the stiffness matrices of $\mathcal{T},\mathcal{T}'$
respectively.  Let $B_S$ be the Schur complement of $B$ with respect to $A$.

If $\mathcal{T}'$ is a fretsaw extension of $\mathcal{T}$, 
then $\lambda_{min}(A,B_S) \geq 1/2$
\end{lemmaE}
\begin{proof}
Suppose that 
$\mathcal{T}=\la m,\{\vv_{i} \}_{i=1}^{m},E',\gamma'\ra$
is the $(\rho,\pi)$-fretsaw extension of
$\mathcal{T}=\la n,\{\vv_{i} \}_{i=1}^{n},E,\gamma\ra$.

Define $M$ to be the $2n\times 2m$ matrix that for all
$i\in [n],j\in [m]$ satisfies
\[
\lmx M_{2i-1} & M_{2i} \\ M_{2j-1} & M_{2j} \rmx =
\begin{cases}
I & i=\pi (j)
\\
0 & \text{otherwise}
\end{cases}
\]
and note that for an element $(i,j)\in E'$, 
$A_{(\pi (i),\pi (j))}=MA_{(i,j)}M^{T}$.  However any element in $\mathcal{T}$
is part of at most two faces, and so can have at most two copies in
$\mathcal{T}'$.  Thus, $A\cleq MBM^{T}\cleq 2A$.

Recalling that $\pi (i)=i$, we note also that $M$ takes the form
$\lmx I & M_{1} \rmx$ for some $(2n-2m)\times 2m$ matrix $M_{1}$.  
Thus, for any $\xx$, we have:
\begin{eqnarray*}
2\xx^T A \xx 
&\geq &
\xx^T M B M^{T} \xx \\
&=&
\lmx \xx^T & M_{1}\xx^T \rmx B \lmx \xx \\ M_{1}^{T}\xx \rmx \\
&\geq&
\xx^TB_S\xx \\
\end{eqnarray*}
where the last inequality holds by Lemma \ref{schurprop}.
\end{proof}

\section{Augmented Spanning Tree}\label{congapp}

We give a proof of Lemma \ref{lowcong}, which we restate
here for convenience.  It is a generalization of an algorithm from \cite{SpielmanTengLinsolve}.

\newtheorem*{lemmacong}{Lemma \ref{lowcong}}
\begin{lemmacong}
There exists an algorithm $S=\mathtt{LowCongestAugment}(G,T,Z,\psi,k)$
that takes a planar graph $G=(V,E)$, a spanning tree $T$ of $G$,
a set $Z$ of pairs of vertices in $V$, an embedding 
$\psi: Z\rightarrow\mathcal{P}(E)$, and an integer $k$.
The algorithm runs in time $O(|E|\log |V| + cong(\psi)|E|)$
and returns a set of edges $S\subseteq E$ of
size at most $k$,
such that there exists an embedding $\pi:Z\rightarrow\mathcal{P}(T\cup S)$
with congestion $O(\frac{1}{k}str(T)cong(\psi))$.
\end{lemmacong}

We need to make use of the following tree decomposition algorithm
from \cite{SpielmanTengLinsolve}:

\begin{theorem}[\cite{SpielmanTengLinsolve} Theorem 8.3]\label{decompose}
There exists a linear-time algorithm
$$
((W_1,...,W_c),\rho) = \mathtt{decompose}(T,E,\eta,k)
$$
that on input
\begin{itemize}
\item a tree $T$ on vertices $V$
\item a set $E$ of edges forming a planar graph on $V$
\item a function $\eta:E\rightarrow \Reals+$
\item a positive integer $k \leq \sum_{e\in E} \eta(e)$
\end{itemize}
outputs sets $W_1,...,W_c\subseteq V$, where $c\leq k$, and a function 
$\rho$ that maps each
edge in $E$ to either a set or pair of sets in $\{W_1,...,W_c\}$ such that:
\begin{itemize}
\item $V = \cup_{i=1}^c W_i$, and for all $i \ne j$, $|W_i\cup W_j|\leq 1$
\item for all $i$, the graph induced by $T$ on $W_i$ (which we denote $T_i$)
  is connected
\item for each $(u,v)\in E$, there are $i,j$ (possibly equal)
  such that $u\in W_i$ and $v\in W_j$ and $\rho(u,v) = \{W_i,W_j\}$
\item the graph 
  $\la [c],\{(i,j): \exists e\in E \text{ s.t. } \rho(e)=\{W_i,W_j\}\} \ra$ 
  is planar
\item for all $W_i$ such that $|W_i|>1$,
$$
\sum_{e\in E:W_i\in\rho(e)} \eta(e) \leq \frac{4}{k}\sum_{e\in E} \eta(e)
$$
\end{itemize}
\end{theorem}

We can now prove the lemma.

\begin{proof}
Here is the {\tt LowCongestAugment} algorithm:

\mybox{
$S=\mathtt{LowCongestAugment}(\la V,E \ra,T,Z,\psi,k)$
\begin{enumerate}
\item For each $e\in E$, 
  define $\eta(e) = \sum_{z:e\in\psi(z)}\sum_{e'\in \psi(z)}|T(e')|$
\item Set $((W_1,...,W_c),\rho) = 
  \mathtt{decompose}(T,E,\eta,\lfloor\frac{k}{3}\rfloor)$.
\item Let $E_{i,j}=\{e\in E: \rho(e)=\{W_i,W_j\}\}$.
  For all nonempty $E_{i,j}$ define
  $$
  s_{i,j} = \arg\min_{e\in E_{i,j}} |T(e)|
  $$

  and let $S$ be the set of all the $s_{i,j}$.
\end{enumerate}
}

Let us analyze the running time.  In step 1, we must compute $\eta(e)$ for all
$e\in E$.  First we compute and record 
$|T(e)|$ for all $e\in E$.  \cite{SpielmanTengLinsolve} gives a method to
do this in time $O(|E|\log |V|)$.  We can then compute each $\eta(e)$
by summing $\sum_{z:e\in\psi(z)}|\psi(z)| \leq cong(\psi)$ of the $|T(e)|$
values.
This gives a total time of $O(|E|\log |V| + cong(\psi)|E|)$ for step 1, 
and the remaining steps clearly run faster than $O(|E|\log |V|)$.

Let also note that $|S|\leq k$.
This follows from the fact that the graph
$\la[c],\{(i,j):E_{i,j}\ne\emptyset\}\ra$ is planar
and has $|S|$ edges.
Since a planar
graph on $c\leq k/3$ vertices cannot have more than 
$3c-6 < k$ edges, we have that 
$|S|\leq k$.

Now let us demonstrate the existence of a low congestion embedding.
For each $(v,w)\in E$, let us define a path $\pi(v,w)$ in $T\cup S$
from $v$ to $w$, as follows:
\begin{itemize}
\item 
If $\rho(v,w)$ is a singleton $\{W_i\}$, then we simply define
$\pi(v,w)=T(v,w)$

Note that $\pi(v,w)\subseteq T_i$.

\item
If $\rho(v,w)$ is a pair $\{W_i,W_j\}$, 
then let $(v',w')=s_{i,j}\in S$
and define $\pi(v,w)=T(v,v')\cup\{(v',w')\}\cup T(w',w)$.

Note that
$$
\pi(v,w)\subseteq T_i\cup T_j\cup \{s_{i,j}\}
$$
and that
$$
|\pi(v,w)| \leq |T(v,w)|+|T(v',w')|+1 \leq 3|T(v,w)|
$$
\end{itemize}

Furthermore for $(v,w) \in S$, 
define $\pi(v,w) = \cup_{e\in \psi(v,w)}\pi(e)$.

Fix an $e_0\in T\cup S$.  By construction of $\pi$, $e$ can be only in a path
$\pi(z)$ if it is either in $S$ or in some subtree $T_k$.

With this in mind, 
define $i_0$ such that $T_{i_0}$ is the
subtree containing $e$.  (There is at most one such tree, but if there is none
then choose $i_0$ arbitrarily.)
Define $i_1,i_2$ such that if $e\in S$ then $e=z_{i_1,i_2}$.
(If $e\not\in S$ then choose $i_1,i_2$ arbitrarily).  
Then for any $z\in Z$, if $e\in \pi(z)$ then 
at least one of the following must hold:
\begin{itemize}
\item For some $e\in\psi(z)$, $W_{i_0} \in \rho(e)$.
\item For some $e\in\psi(z)$, $\rho(e)=\{W_{i_1},W_{i_2}\}$.  (In particular,
$W_{i_1}\in\rho(e)$.)
\end{itemize}

Thus we can bound the congestion of $\pi$ on $e_0$:
\begin{eqnarray*}
\sum_{z:e_0\in\pi(z)}|\pi(z)| 
&\leq&
\sum_{z:\{e\in\psi(z):W_{i_0}\in\rho(e)\}\ne\emptyset} |\pi(z)| 
\sum_{z:\{e\in\psi(z):W_{i_1}\in\rho(e)\}\ne\emptyset} |\pi(z)| 
\end{eqnarray*}

We note that for any $i$:
\begin{eqnarray*}
\sum_{z:\{e\in\psi(z):W_i\in\rho(e)\}\ne\emptyset} |\pi(z)| 
&\leq&
\sum_{e:W_i\in\rho(e)}\sum_{z:e\in\psi(z)} |\pi(z)| \\
&\leq&
\sum_{e:W_i\in\rho(e)}\sum_{z:e\in\psi(z)}\sum_{e'\in \psi(z)}|\pi(e')| \\ 
&\leq&
3\sum_{e:W_i\in\rho(e)}\sum_{z:e\in\psi(z)}\sum_{e'\in \psi(z)}|T(e')| \\ 
&=&
3\sum_{e:W_i\in\rho(e)}\eta(e) \\
&\leq&
\frac{12}{k}\sum_{e\in E}\eta(e) \\
\end{eqnarray*}

Thus:
\begin{eqnarray*}
\sum_{z:e_0\in\pi(z)}|\pi(z)| 
&\leq&
\frac{24}{k}\sum_{e\in E}\eta(e) \\
&\leq&
\frac{24}{k}\sum_{e\in E}\sum_{z:e\in\psi(z)}\sum_{e'\in \psi(z)}|T(e')| \\ 
&=&
\frac{24}{k}\sum_{e'\in E}|T(e')|\sum_{z:e'\in\psi(z)}|\psi(z)| \\ 
&\leq&
\frac{24}{k}\sum_{e'\in E}|T(e')|cong(\psi) \\ 
&=&
\frac{24}{k}str(T)cong(\psi) \\ 
\end{eqnarray*}
\end{proof}

Since the choice of $e_0$ was arbitrary, we have
$$
cong(\pi) \leq \frac{24}{k}str(T)cong(\psi)
$$

\end{document}